*To my daughters Vianne and Natalie*

# AN ALGEBRAIC CONSTRUCTION OF HYPERBOLIC PLANES OVER A EUCLIDEAN ORDERED FIELD

**Nicholas Phat Nguyen**[1]

**Abstract.** Using concepts and techniques of bilinear algebra, we construct hyperbolic planes over a euclidean ordered field that satisfy all the Hilbert axioms of incidence, order and congruence for a basic plane geometry, but for which the hyperbolic version of the parallel axiom holds rather than the classical Euclidean parallel postulate.

**1. INTRODUCTION.** Most readers are probably familiar with the Poincare model of a hyperbolic plane, based on either the upper half-plane or the interior of the unit disc in the standard Euclidian plane (which can also be thought of as the space of complex numbers). Such a hyperbolic plane is a 2-dimensional real manifold (or a 1-dimensional complex manifold) where a suitable concept of lines could be defined so that the hyperbolic version of the Euclidean parallel axiom holds: for any point not on a given line, there are two or more lines (in fact an infinite number of lines) passing through that point and not intersecting with the given line.

We can construct a hyperbolic plane using a similar Poincare model applied to an affine plane over an euclidean ordered field. See [1] at chapter 7. Such a construction gives us a model that satisfies all the basic axioms of incidence, order and congruence for lines, segments and angles for a basic Hilbert plane geometry, but for which the hyperbolic version of the parallel axiom holds rather than the classical Euclidean parallel postulate.

[1] Email: nicholas.pn@gmail.com

In the following, we will outline an alternative construction of such a Poincare model over any euclidean ordered field using methods from bilinear algebra (i.e. the study of symmetric bilinear forms or quadratic forms over a field). Such an algebraic construction provides new methods that could help in the exploration of general hyperbolic geometry. The construction also does not rely on any visual verification for its definitions or its proofs, and gives us more precise information, such as a natural parametrization of the set of lines through a point that are not parallel to a given line.

Recall that an ordered field is a field endowed with a total ordering relation compatible with the additive and multiplicative operations of the field. An ordered field could be defined by giving a total ordering $\leq$ that satisfies the following basic properties associated with the usual ordering on the real line:

- if $a \leq b$ then $a + c \leq b + c$ *for any* $c$
- if $0 \leq a$ and $0 \leq b$ then $0 \leq a\, b$

An ordered field is said to be euclidean if every positive number is a square. Examples of euclidean fields include the field of real numbers **R**, the field of all algebraic numbers contained in **R**, and the field of all real constructive numbers (real numbers that can be constructed from the rational numbers using ruler and compass constructions). Given any ordered field, there are always ordered extensions of that field which are euclidian. So there are infinitely many euclidean ordered fields beyond the subfields of **R** cited in our examples above.

Let $K$ be an ordered field, and consider the n-dimensional vector space $K^n$ of all $n$-tuples with coefficients in $K$. The standard dot product $\boldsymbol{x.y} = x_1y_1 + \ldots + x_ny_n$ is a symmetric bilinear form on $K^n$. This symmetric bilinear form is non-degenerate, in fact anisotropic in the sense that $\boldsymbol{x.x}$ is never zero unless the vector $\boldsymbol{x}$ itself is zero. The fact that the standard dot product is anisotropic for any $K^n$ is a characteristic property of ordered field, as discovered by Emil Artin and Otto Schreier (a field can be ordered if and only if it is formally real).

In the following, we will start out with a vector space $E$ of dimension 2 over a euclidean ordered field $K$. We will assume that $E$ has a symmetric bilinear form $(x \mid y)$ that is positive definite, meaning that the product $(x \mid x)$ is always $> 0$ unless $x$ is the zero vector. Such a positive definite form is obviously anisotropic. The standard dot product on $K^2$ is a prime example. Accordingly, we will use the dot product notation for the bilinear form on $E$ to emphasize the analogy.

Before we begin, we note here a potential difference in terminology. A regular quadratic space of dimension 2 is either anisotropic or isometric to a standard space known as an Artinian or hyperbolic plane which has two linearly independent lines of isotropic vectors. Relative to a basis formed of the two linearly independent isotropic vectors, the matrix of the bilinear form is a 2 by 2 symmetric matrix that has zeros in the diagonal and a non-zero number in the cross diagonal. The literature on bilinear algebra more commonly refers to such a space as a hyperbolic plane, but in order to avoid confusion, we will refer to such a space as an Artinian plane, which is another name that is sometimes used for it. The late geometer Marcel Berger championed the name "Artinian plane" instead of "hyperbolic plane." He regarded the term "hyperbolic plane" often used in bilinear algebra as undesirable, because the term hyperbolic plane can of course also mean a hyperbolic manifold of dimension 2 or a space that satisfies the axioms of hyperbolic plane geometry. See [5] at paragraph 13.1.4.4.

**2. THE SPACE OF CYCLES AND THE CYCLE PAIRING PRODUCT.** Consider the set $\mathcal{F}$ of all functions $p$ from $E$ to $K$ of the form $p(\boldsymbol{X}) = a\boldsymbol{X}.\boldsymbol{X} + \boldsymbol{b}.\boldsymbol{X} + c$, where $\boldsymbol{X}$ and $\boldsymbol{b}$ are vectors in $E$, and $a$ and $c$ are elements in the field $K$.[2] The set $\mathcal{F}$ can naturally be endowed with the

---

[2] If we regard the coefficients of the vector $\boldsymbol{X}$ with respect to a fixed basis of $E$ as variables, the function $p$ is at most a second-degree polynomial in 2 variables. Because the ordered field $K$ has

structure of a $K$-vector space of dimension $2 + 2 = 4$, being parametrized by the vector $\boldsymbol{b}$ and the elements $a$ and $c$. We will refer to such a function $p$ (if it is not the constant zero function) as a 2-cycle, 1-cycle, or 0-cycle depending on whether the degree of $p$ is 2 (coefficient $a$ is nonzero), 1 ($a$ is zero but $\boldsymbol{b}$ is nonzero, or 0 (both $a$ and $\boldsymbol{b}$ are zero).

Aside from the natural structure of a vector space over $K$, we can also endow $\mathcal{F}$ with a symmetric bilinear form $<_,_>$ as follows. Given $p = a\boldsymbol{X}.\boldsymbol{X} + \boldsymbol{b}.\boldsymbol{X} + c$ and $p^* = a^*\boldsymbol{X}.\boldsymbol{X} + \boldsymbol{b^*}.\boldsymbol{X} + c^*$, we define $<p,p^*>$ as $\boldsymbol{b}.\boldsymbol{b^*} - 2ac^* - 2a^*c$.

This scalar product is clearly symmetric and bilinear. Moreover, it is non-degenerate, because $\mathcal{F}$ with this scalar product is isometric to the orthogonal sum of $E$ and an Artinian plane. We will refer to this fundamental scalar product on $\mathcal{F}$ as the cycle pairing or cycle product.

The reader can easily verify the following properties from the formula for the cycle pairing product.

- All 0-cycles are isotropic.
- No 1-cycle is isotropic.
- A 2-cycle $p(\boldsymbol{X}) = a\boldsymbol{X}.\boldsymbol{X} + \boldsymbol{b}.\boldsymbol{X} + c$ is isotropic if and only if $\boldsymbol{b}.\boldsymbol{b} - 4ac = 0$, or equivalently, if $p(\boldsymbol{X}) = a(\boldsymbol{X} + b/2a).(\boldsymbol{X} + b/2a)$.

We will informally refer to such an isotropic 2-cycle as a zero circle centered at the point $(-\boldsymbol{b}/2a)$, in analogy with the familiar circle equation in the Euclidean plane $\mathbf{R}^2$. For each point $\boldsymbol{u}$ in the plane $E$, we write $q(\boldsymbol{u}) = \boldsymbol{X}.\boldsymbol{X} - 2\boldsymbol{u}.\boldsymbol{X} + \boldsymbol{u}.\boldsymbol{u}$ for the normalized zero circle centered at $\boldsymbol{u}$. Note that $\boldsymbol{u}$ is a point or vector in $E$, but $q(\boldsymbol{u})$ is a 2-cycle element of $F$.

Now take a simple 1-cycle $p$ of the form $p(\boldsymbol{X}) = \boldsymbol{i}.\boldsymbol{X}$ where $\boldsymbol{i}$ is a nonzero vector in $E$. Note that the norm $<p, p>$ of the cycle $p$ is equal to $\boldsymbol{i}.\boldsymbol{i}$ and therefore strictly positive. To characteristic 0 and therefore an infinite number of elements, such a function $p$ is a zero function if and only if $a$, $\boldsymbol{b}$, and $c$ are all zero.

simplify notation, we will assume that $p$ has unit norm, which is possible by scaling because the field $K$ is euclidean.

We define $\mathcal{H}$ to be the subset of $E$ consisting of all points $\boldsymbol{u}$ in $E$ such that $<p, q(\boldsymbol{u})> = -2\boldsymbol{i}.\boldsymbol{u}$ is $> 0$. That set corresponds to an open half plane in $E$ bounded by the line with equation $\boldsymbol{i}.\boldsymbol{X} = 0$.

We now defined the hyperbolic lines in the set $\mathcal{H}$ as follows. Any such line is a nonempty subset of $\mathcal{H}$ consisting of all points $\boldsymbol{u}$ in $\mathcal{H}$ such that $<m, q(\boldsymbol{u})> = 0$ for some nonisotropic cycle $m$ orthogonal to $p$, that is to say $<m, p> = 0$. We will informally refer to such points $\boldsymbol{u}$ as the zero points of $m$ in $\mathcal{H}$.

A nonisotropic cycle $m$ orthogonal to $p$ may not have any zero point in $E$, but if $m$ has a zero point in $E$, then it also has an infinite number of zero points in $\mathcal{H}$.

**<u>Proposition 1</u>:** *If a nonisotropic cycle m orthogonal to p has a zero point in E, then m has positive norm. Moreover, m has an infinite number of zero points in $\mathcal{H}$.*

*Proof.* With a suitable base of orthogonal vectors, the matrix of the cycle pairing has the diagonal form [1, 1, 1, –1]. A nonisotropic cycle of negative norm would give us a one-dimensional subspace isomorphic to the one-dimensional quadratic subspace [–1], whose orthogonal complement in $\mathcal{F}$ would be an anisotropic space of dimension 3 isomorphic to the positive definite space [1, 1, 1]. Because that complement is anisotropic, it contains no isotropic 2-cycle, and therefore there is no point $\boldsymbol{u}$ such that $q(\boldsymbol{u})$ is orthogonal to the given cycle.

So if a cycle $m$ has a zero point in $\mathcal{H}$, it must have positive norm. All the different isotropic lines in the orthogonal complement of $m$ correspond to points on the projective conic $U^2 + V^2 - W^2 = 0$, and there are an infinite number of these points.[3] There are exactly two isotropic points orthogonal to both $p$ and $m$ because the vector subspace generated by

[3] The points of this projective conic can be put in bijective correspondence with the points on a projective line.

$p$ and $m$ is anisotropic and so its orthogonal complement must be an Artinian plane. If $\boldsymbol{u}$ is a zero point of $m$ in $E$ not orthogonal to $p$, then the cycle $r = q(\boldsymbol{u}) - 2 <p, q(\boldsymbol{u})>p$ is a normalized zero circle whose cycle product with $p$ has the opposite sign and whose cycle product with $m$ is also zero. Therefore, the zero points of $m$ in $E$ not orthogonal to $p$ can be divided into two subsets of equal size. That means the hyperbolic line defined by $m$ must have an infinite number of points. ■

Note that if a nonisotropic cycle $m$ defines a line in $\mathcal{H}$, then any multiple $\mu m$ of $m$ ($\mu \neq 0$) defines the same line. The complement of such a line in $\mathcal{H}$ can naturally be divided into two sides, each side consisting all points $\boldsymbol{u}$ such that $<m, q(\boldsymbol{u})>$ has the same sign.

## 3. THE HYPERBOLIC PLANE – INCIDENCE AXIOMS.

We now show that the set $\mathcal{H}$ with the lines as defined above satisfies the Hilbert axioms for a hyperbolic plane. We begin with the incidence axioms.

*(Incidence Axiom 1) For every pair of distinct points* $\boldsymbol{u}$ *and* $\boldsymbol{v}$*, there is a unique line containing both* $\boldsymbol{u}$ *and* $\boldsymbol{v}$*.*

Let $\boldsymbol{u}$ and $\boldsymbol{v}$ be two distinct points in $\mathcal{H}$. We claim that the cycles $p$, $q(\boldsymbol{u})$ and $q(\boldsymbol{v})$ are linearly independent. Because $q(\boldsymbol{u})$ and $q(\boldsymbol{v})$ are clearly not proportional, linear dependence in this case would mean there is an equation $p = aq(\boldsymbol{u}) + bq(\boldsymbol{v})$ for some nonzero numbers $a$ and $b$ in the field $K$. Because $p$ is a 1-cycle, this means $a + b = 0$, or that these two numbers have opposite signs. Let's say that $a$ is $> 0$. If we take the cycle product of both sides with $q(\boldsymbol{v})$, we see that $<p, q(\boldsymbol{v})> = a<q(\boldsymbol{u}), q(\boldsymbol{v})>$. The left hand side is $> 0$ by definition of $\mathcal{H}$. For the right hand side, note that $<q(\boldsymbol{u}), q(\boldsymbol{v})> = 4\boldsymbol{u}.\boldsymbol{v} - 2\boldsymbol{u}.\boldsymbol{u} - 2\boldsymbol{v}.\boldsymbol{v} = -2(\boldsymbol{u} - \boldsymbol{v}).(\boldsymbol{u} - \boldsymbol{v})$ is $< 0$ by hypothesis. So the right hand side $a<q(\boldsymbol{u}), q(\boldsymbol{v})>$ is $< 0$. This would be a contradiction, so the cycles $p$, $q(\boldsymbol{u})$ and $q(\boldsymbol{v})$ must be linearly independent.

The subspace of $\mathcal{F}$ generated by the cycles $p$, $q(\boldsymbol{u})$ and $q(\boldsymbol{v})$ therefore has dimension 3. Moreover we can readily check that the cycle pairing on this subspace has discriminant $<$ 0, which means the cycle pairing on this subspace is regular. Therefore the orthogonal complement of this subspace in $\mathcal{F}$ is a one-dimensional subspace generated by a nonisotropic cycle $m$. By construction $m$ is orthogonal to $p$, $q(\mathbf{u})$ and $q(\mathbf{v})$, and therefore gives us a line in $\mathcal{H}$ that contains $\mathbf{u}$ and $\mathbf{v}$. This line is unique, since any such line by definition must correspond to a cycle that is orthogonal to $p$, $q(\mathbf{u})$ and $q(\mathbf{v})$, and therefore that cycle must be proportional to the cycle $m$. Proportional cycles obviously define the same line. ■

*(Incidence Axiom 2) Every line contains at least two points.*

Proposition 1 tells us that every line in fact contains an infinite number of points. ■

*(Incidence Axiom 3) There are at least three points that are not collinear.*

We can chose a basis $\mathbf{s}$ and $\mathbf{t}$ in $E$ such that $\mathbf{i}.\mathbf{s}$ and $\mathbf{i}.\mathbf{t}$ are $< 0$. For any point $\mathbf{u}$ in $\mathcal{H}$, the points $\mathbf{u} + \mathbf{s}$, $\mathbf{u} + \mathbf{t}$ are also in $\mathcal{H}$. Moreover, the points $\mathbf{u}$, $\mathbf{u} + \mathbf{s}$, and $\mathbf{u} + \mathbf{t}$ are not all collinear. Otherwise we would have an nonisotropic cycle $m$ that is orthogonal to $p$, $q(\mathbf{u})$, $q(\mathbf{u} + \mathbf{s})$ and $q(\mathbf{u} + \mathbf{t})$. But we can check by straight-forward computations that these four cycles generate the entire space $\mathcal{F}$ of all cycles. The cycle pairing is nondegenerate, so the only vector in $\mathcal{F}$ cycle orthogonal to all these four cycles is the zero vector. ■

**4. THE HYPERBOLIC PLANE – ORDER AXIOMS.** We need to define what it means for a point to lie between two other points on the same line.

Let $\mathbf{u}$, $\mathbf{v}$ and $\mathbf{w}$ be three distinct collinear points in $\mathcal{H}$. The fact that these three points are collinear means first of all that the four cycles $p$, $q(\mathbf{u})$, $q(\mathbf{v})$ and $q(\mathbf{w})$ are not linearly independent. Otherwise they would generate the entire space $\mathcal{F}$ of cocycles, and in that case there would be no nonisotropic cycle orthogonal to all of them.

We claim that $q(\mathbf{u})$, $q(\mathbf{v})$ and $q(\mathbf{w})$ must be linearly independent, because otherwise we would have an equation $aq(\mathbf{u}) + bq(\mathbf{v}) + cq(\mathbf{w}) = 0$ with $abc \neq 0$. Two of the three coefficient numbers must have the same sign. Let's say $a$ and $b$ have the same sign. If we now take cycle product of both sides of the equation with $q(\mathbf{w})$, the right hand side is obviously zero. On the left hand side we have $a<q(\mathbf{u}), q(\mathbf{w})> + b<q(\mathbf{v}), q(\mathbf{w})> \neq 0$ because it is the sum of two numbers of the same sign. (Recall that for any two different points $\mathbf{x}$ and $\mathbf{y}$ in $E$, $<q(\mathbf{x}), q(\mathbf{y})> = -2(\mathbf{x} - \mathbf{y}).(\mathbf{x} - \mathbf{y})$ is $< 0$.)

Accordingly, the linear dependence of the four cycles $p$, $q(\boldsymbol{u})$, $q(\boldsymbol{v})$ and $q(\boldsymbol{w})$ implies that we must have an equation $p = aq(\boldsymbol{u}) + bq(\boldsymbol{v}) + cq(\boldsymbol{w}) = 0$ with $abc \neq 0$. We claim that the three coefficient numbers cannot all have the same sign. Otherwise, by taking cycle product of both sides with $p$, we deduce that all three numbers must be positive. On the other hand, by taking the cycle product of both sides with $q(\boldsymbol{w})$, we see that the left hand side is strictly positive (by definition of the set $\boldsymbol{\mathcal{H}}$) while the right hand side is strictly negative because $<q(\boldsymbol{u}), q(\boldsymbol{w})>$ and $<q(\boldsymbol{v}), q(\boldsymbol{w})>$ are both $< 0$.

Therefore exactly two of the three coefficients have the same sign. The remaining coefficient has a different sign, and we say that the point corresponding to that coefficient lies between the other two points.

Our definition immediately implies the following two axioms:

*(Order Axiom 1) For three collinear points **u**, **v** and **w**, if **w** is between **u** and **v**, then **w** is also between **v** and **u**.*

*(Order Axiom 2) For any three collinear points **u**, **v** and **w**, exactly one point is between the other two points.*

We now consider the following more difficult axioms.

*(Order Axiom 3) For any two points **s** and **t**, we can find at least three points **u**, **v** and **w** on the line passing through **s** and **t** such that **s** is between **u** and **t**, **v** is between **s** and **t**, and **t** is between **s** and **w**.*

We will later prove that for any two lines $L$ and $M$ in $\boldsymbol{\mathcal{H}}$, there is a transformation of that maps the line $L$ bijectively onto the line $M$ and respecting the betweenness relationship of any three points. Assuming this, it is enough for us to show that Order Axiom 3 applies to one specific line.

Consider the line defined by the 1-cycle $\boldsymbol{j}.\boldsymbol{X}$, where $\boldsymbol{j}$ is a nonzero vector in $\boldsymbol{E}$ such that $\boldsymbol{i}.\boldsymbol{j} = 0$. This line consists of all points in $E$ of the form $-2t\boldsymbol{i}$ where $t$ runs through all numbers $> 0$. On this line, there is a natural betweenness induced by the given ordering of $K$. That natural betweenness clearly satisfies Order Axiom 3. We will show that the natural betweenness relationship on this particular line is the same as the betweenness relationship defined above.

Let $\boldsymbol{u}$, $\boldsymbol{v}$ and $\boldsymbol{w}$ be three distinct points on this line. Under our definition, to say that $\boldsymbol{v}$ lies between $\boldsymbol{u}$ and $\boldsymbol{w}$ means there is an equation $p = xq(\boldsymbol{u}) + yq(\boldsymbol{v}) + zq(\boldsymbol{w})$ where $x$, $y$ and $z$ are nonzero numbers such that $x$ and $z$ have the same sign. Let $\boldsymbol{u} = -2A\boldsymbol{i}$, $\boldsymbol{v} = -2B\boldsymbol{i}$, and $\boldsymbol{w} = -2C\boldsymbol{i}$. The equation $p = xq(\boldsymbol{u}) + yq(\boldsymbol{v}) + zq(\boldsymbol{w})$ means that the cycle $\boldsymbol{i.X}$ is the same as the cycle $x(\boldsymbol{X.X} - 2\boldsymbol{u.X} + \boldsymbol{u.u}) + y(\boldsymbol{X.X} - 2\boldsymbol{v.X} + \boldsymbol{v.v}) + z(\boldsymbol{X.X} - 2\boldsymbol{w.X} + \boldsymbol{w.w})$. By equating the coefficients of the cycles on both sides, we have the following three linear equations in three unknowns x, y and z:

$$x + y + z = 0$$

$$-2Ax - 2By - 2Cz = 1$$

$$A^2x + B^2y + C^2z = 0$$

From linear algebra, we know that up to a common factor, $x$ and $z$ are proportional to $(C^2 - B^2)$ and $(B^2 - A^2)$. That means the numbers $x$ and $z$ have the same sign if and only if we have $(C^2 > B^2 > A^2)$ or $(C^2 < B^2 < A^2)$. But these conditions on the numbers $A$, $B$ and $C$ mean that $\boldsymbol{v}$ lies between $\boldsymbol{u}$ and $\boldsymbol{w}$ in the natural ordering on that line. ■

*(Order Axiom 4 – Pasch's Axiom) Let* ***u, v, w*** *be three distinct points in* $\mathcal{H}$*, and let M be a line in* $\mathcal{H}$ *that does not pass through any of these three points. If M passes through a point between* ***u*** *and* ***v****, then M also passes through a point either between* ***u*** *and* ***w****, or between* ***v*** *and* ***w.***

Let $m$ be a nonisotropic cycle defining the line $M$. We claim that in order for the line $M$ to pass through a point between $\boldsymbol{r}$ and $\boldsymbol{s}$, it is necessary and sufficient that the numbers $<m, q(\boldsymbol{r})>$ and $<m, q(\boldsymbol{s})>$ have opposite signs. Such a result will immediately establish the Pasch's axiom above. Indeed, if $m$ passes through a point between $\boldsymbol{u}$ and $\boldsymbol{v}$, then that means $< m, q(\boldsymbol{u}) >$ and $< m, q(\boldsymbol{v}) >$ have opposite signs. The number $< m, q(\boldsymbol{w}) >$ obviously would have the same sign with either $< m, q(\boldsymbol{u}) >$ or $< m, q(\boldsymbol{v}) >$ but not both. Let us say that $< m, q(\boldsymbol{w}) >$ has the same sign with $< m, q(\boldsymbol{v}) >$ and the opposite sign with $< m, q(\boldsymbol{u}) >$. Then the line $M$ will pass through a point between $\boldsymbol{u}$ and $\boldsymbol{w}$, but not through any point between $\boldsymbol{v}$ and $\boldsymbol{w}$.

To prove our claim, assume first that the line $M$ intersects the line $L$ that goes through $\boldsymbol{r}$ and $\boldsymbol{s}$. We can show that the intersection point $\boldsymbol{t}$ lies between $\boldsymbol{r}$ and $\boldsymbol{s}$ if and only if the numbers $< m, q(\boldsymbol{r}) >$ and $< m, q(\boldsymbol{s}) >$ have opposite signs. Specifically, we know that the four cycles $p$, $q(\boldsymbol{r})$, $q(\boldsymbol{s})$ and $q(\boldsymbol{t})$ are linearly dependent and that we must have an equation $p = aq(\boldsymbol{r}) + bq(\boldsymbol{s}) + cq(\boldsymbol{t}) = 0$ with $abc \neq 0$.

Now take the cycle pairing with $m$. On the left hand side we have $< m, p > = 0$ because by definition any cycle defining a line is orthogonal to $p$. On the right hand side, we have the sum $a< m, q(\boldsymbol{r}) > + b< m, q(\boldsymbol{s}) >$ because $< m, q(\boldsymbol{t}) > = 0$ (since the line $M$ passes through $\boldsymbol{t}$). Therefore $a< m, q(\boldsymbol{r}) > + b< m, q(\boldsymbol{s}) > = 0$, and for the numbers $a$ and $b$ to have the same sign (that is to say for $\boldsymbol{t}$ to be between $\boldsymbol{r}$ and $\boldsymbol{s}$), it is necessary and sufficient that the numbers $< m, q(\boldsymbol{r}) >$ and $< m, q(\boldsymbol{s}) >$ have opposite signs.

We now show that if the numbers $< m, q(\boldsymbol{r}) >$ and $< m, q(\boldsymbol{s}) >$ have opposite signs, then the line $M$ must intersect the line $L$ through $\boldsymbol{r}$ and $\boldsymbol{s}$. Let $\ell$ be a nonisotropic cycle that defines the line $L$. We claim we can find numbers $x$, $y$ and $z$ such that the linear combination $xp + yq(\boldsymbol{r}) + zq(\boldsymbol{s})$ is cycle orthogonal to both $\ell$ and $m$, and such that $y$ and $z$ are nonzero numbers of the same sign.

Such a linear combination is of course orthogonal to $\ell$ by definition of the line $L$. Because $m$ is orthogonal to $p$ by definition, we just need to make sure that $m$ is orthogonal to $yq(\boldsymbol{r}) + zq(\boldsymbol{s})$. Note that by the hypothesis of $< m, q(\boldsymbol{r}) >$ and $< m, q(\boldsymbol{s}) >$ having opposite signs, we can find nonzero numbers $y$ and $z$ of the same sign such that we have

$$< m, yq(\boldsymbol{r}) + zq(\boldsymbol{s}) > = y< m, q(\boldsymbol{r}) > + z< m, q(\boldsymbol{s}) > = 0.$$

Assume $y$ and $z$ chosen, we can then find a number $x$ such that the linear combination $xp + yq(\boldsymbol{r}) + zq(\boldsymbol{s})$, which is a 2-cycle, has zero norm. That zero norm condition is equivalent to the equation $\langle xp + yq(\boldsymbol{r}) + zq(\boldsymbol{s}), xp + yq(\boldsymbol{r}) + zq(\boldsymbol{s}) \rangle = 0$, or

$$x^2 + 2x\langle p, yq(\boldsymbol{r}) + zq(\boldsymbol{s})\rangle + 2yz\langle q(\boldsymbol{r}), q(\boldsymbol{s})\rangle = 0.$$

Such a monic quadratic equation must have a root $x$ in the euclidean field $K$ because the constant term $2yz\langle q(\boldsymbol{r}), q(\boldsymbol{s})\rangle$ is strictly negative. Indeed, the number $yz$ is $> 0$ because $y$ and $z$ are nonzero numbers of the same sign, and we know from earlier computation that $\langle q(\boldsymbol{r}), q(\boldsymbol{s})\rangle$ is a strictly negative number for any two points $\boldsymbol{r}$ and $\boldsymbol{s}$ in $\mathcal{H}$.

The 2-cycle $xp + yq(\boldsymbol{r}) + zq(\boldsymbol{s})$ therefore is isotropic and must be equal to $\alpha q(\boldsymbol{t})$ for some vector $\boldsymbol{t}$ in $E$ and some nonzero coefficient $\alpha$. By construction $q(\boldsymbol{t})$ is cycle orthogonal to both $\ell$ and $m$.

If $\langle p, q(\boldsymbol{t})\rangle$ is strictly positive, then $\boldsymbol{t}$ is a point in $\mathcal{H}$, and we are done. If not, then the reflected 2-cycle $q(\boldsymbol{t}) - 2\langle p, q(\boldsymbol{t})\rangle p$ is a normalized zero circle that represents a point in $\mathcal{H}$ lying on both $M$ and $L$, and also lying between $\boldsymbol{r}$ and $\boldsymbol{s}$. ■

Given a point $\boldsymbol{u}$ on a line $L$, we can divide the complement of $\boldsymbol{u}$ in $L$ into two sides as follows. Let $\ell$ be a nonisotropic cycle defining the line $L$. The cycles $p$, $q(\boldsymbol{u})$ and $\ell$ generate a regular subspace of dimension 3 in $\mathcal{F}$. The orthogonal complement of that subspace is a one-dimensional subspace generated by a nonisotropic cycle $n$. For any point $\boldsymbol{v} \neq \boldsymbol{u}$ on $L$, we must have $\langle n, q(\boldsymbol{v})\rangle \neq 0$, because otherwise the cycle $n$ would be orthogonal to $p$, $q(\boldsymbol{u})$, $q(\boldsymbol{v})$ and $\ell$ and hence orthogonal to the whole space $\mathcal{F}$. The complement of $\boldsymbol{u}$ in $L$ can then be divided into two sides according to whether $\langle n, q(\boldsymbol{v})\rangle$ is $> 0$ or $< 0$.

If $\boldsymbol{w}$ is any point $\neq \boldsymbol{u}$ on $L$, we will refer to the set containing $\boldsymbol{u}$ and the side containing $\boldsymbol{w}$ as the ray $[\boldsymbol{uw}, \infty)$.

**5. THE HYPERBOLIC PLANE – CONGRUENCE AXIOMS.** We want to construct a group of transformations of $\mathcal{H}$ that map lines to lines and preserve the betweenness relationship. Once we have such a group, the definition of congruence follows naturally.

Consider all the orthogonal transformations of the vector space $\mathcal{F}$. These are invertible linear transformations of $\mathcal{F}$ that preserve the cycle pairing. Let's focus on the orthogonal transformations that fix the cycle $p$. If $T$ is such a transformation, then for any point $\boldsymbol{u}$ in $\mathcal{H}$, $T(q(\boldsymbol{u}))$ is an isotropic cycle whose cycle product with $p$ is $> 0$. $T(q(\boldsymbol{u}))$ is not a 0-cycle because any 0-cycle is orthogonal to $p$. So $T(q(\boldsymbol{u}))$ must be an isotropic 2-cycle, and therefore can be written as $\lambda q(\boldsymbol{v})$ for some point $\boldsymbol{v}$ in $E$. For $\boldsymbol{v}$ to be in $\mathcal{H}$, $<p, q(\boldsymbol{v})>$ must be $> 0$, and hence it is necessary and sufficient that $\lambda > 0$. We will call any orthogonal transformation $T$ a proper transformation if it fixes $p$ and has the property that for any $\boldsymbol{u}$ in $\mathcal{H}$, we have $T(q(\boldsymbol{u})) = \lambda q(\boldsymbol{v})$ for some $\lambda > 0$. Such a proper transformation clearly induces a bijection of   to itself, and gives us a transformation of $\mathcal{H}$ that maps lines to lines and preserves the betweenness relationship.

Any orthogonal transformation of $\mathcal{F}$ can be written as a product of a finite number of basic transformations called reflections. (This is known as the Cartan – Dieudonne theorem.) Each nonisotropic cycle of $\mathcal{F}$ defines such a reflection as follows. Let $\boldsymbol{x}$ be a nonisotropic cycle of $\mathcal{F}$, then by orthogonal decomposition each cycle $\boldsymbol{v}$ of $\mathcal{F}$ can be expressed uniquely as a sum $\boldsymbol{v} = a\boldsymbol{x} + \boldsymbol{y}$, where $\boldsymbol{y}$ is a cycle in $\mathcal{F}$ orthogonal to $\boldsymbol{x}$. The reflection defined by $\boldsymbol{x}$ maps $\boldsymbol{v}$ to $\boldsymbol{v}' = -a\boldsymbol{x} + \boldsymbol{y}$. In other words, the reflection leaves invariant the orthogonal complement of $\boldsymbol{x}$ and maps $\boldsymbol{x}$ to $-\boldsymbol{x}$. Such a reflection is also called a hyperplane reflection along $\boldsymbol{x}$ or across the orthogonal complement of $\boldsymbol{x}$. Note that all nonzero scalar multiples of the same nonisotropic cycle define the same reflection.

It follows that any transformation $T$ leaving the cycle $p$ invariant can be expressed as a product of reflections defined by nonisotropic cycles that are orthogonal to $p$.

<u>Proposition 2</u>: *Any reflection defined by a nonisotropic cycle that is orthogonal to p and that has positive norm is a proper transformation.*

*Proof.* Let $t = a\boldsymbol{X}.\boldsymbol{X} + \boldsymbol{b}.\boldsymbol{X} + c$ be a cycle orthogonal to $p$ such that its norm $<t, t> = \boldsymbol{b}.\boldsymbol{b} - 4ac$ is $> 0$. To simplify notation, we will assume that $<t, t> = 1$. Because the field $K$ is euclidean, such scaling is possible. Note also that the condition $<t, p> = 0$ means $\boldsymbol{b}.\boldsymbol{i} = 0$

Let $\boldsymbol{u}$ be any point in $\mathcal{H}$. The normalized zero circle centered at $\boldsymbol{u}$ is the 2-cycle $q(\boldsymbol{u}) = \boldsymbol{X.X} - 2\boldsymbol{u.X} + \boldsymbol{u.u}$. We have $<q(\boldsymbol{u}), t> = -2\boldsymbol{u.b} - 2a\boldsymbol{u.u} - 2c$.

The reflection defined by $t$ is the transformation $R(q(\boldsymbol{u})) = q(\boldsymbol{u}) - 2 <q(\boldsymbol{u}), t> t$. It is a proper transformation precisely when the coefficient of the quadratic term $\boldsymbol{X.X}$ in $R(q(\boldsymbol{u}))$ is a strictly positive number. That coefficient is $1 - 2(-2\boldsymbol{u.b} - 2a\boldsymbol{u.u} - 2c)a = 1 + 4a\boldsymbol{u.b} + 4aa\boldsymbol{u.u} + 4ac = 1 + (2a\boldsymbol{u} + \boldsymbol{b}).(2a\boldsymbol{u}, \boldsymbol{b}) - \boldsymbol{b.b} + 4ac = (2a\boldsymbol{u} + \boldsymbol{b}).(2a\boldsymbol{u}, \boldsymbol{b})$.

Clearly $(2a\boldsymbol{u} + \boldsymbol{b}).(2a\boldsymbol{u}, \boldsymbol{b})$ will always be $> 0$ unless $(2a\boldsymbol{u} + \boldsymbol{b})$ is the zero vector in $E$. But in our situation $(2a\boldsymbol{u} + \boldsymbol{b})$ will always be a nonzero vector, as explained below.

- If $a = 0$, then $\boldsymbol{b.b} = 1$ so $(2a\boldsymbol{u} + \boldsymbol{b}) = \boldsymbol{b}$ is a nonzero vector.
- If $a \neq 0$, then note that $(2a\boldsymbol{u} + \boldsymbol{b}).\boldsymbol{i} = 2a\boldsymbol{u.i} + \boldsymbol{b.i} = 2a\boldsymbol{u.i}$ because $t$ is orthogonal to $p$. Recall that $-2\boldsymbol{u.i} > 0$ for any point $\boldsymbol{u}$ in according to our definition of $\mathcal{H}$, so that $2a\boldsymbol{u.i} \neq 0$. Therefore $(2a\boldsymbol{u} + \boldsymbol{b})$ must be a nonzero vector.

Consequently, the reflection defined by $t$ is a proper transformation. ■

Let $G$ be the group of proper orthogonal transformations as defined above. The group $G$ would include reflections defined by nonisotropic cycles that are orthogonal to $p$ and have positive norm (such as the cycles corresponding to lines in $\mathcal{H}$), as well as any composition of such reflections. The set $\mathcal{H}$ is stable under the action of $G$, and hence each transformation in $G$ induces a transformation of $\mathcal{H}$ which we will call a congruence transformation. For convenience, we will often use the same letter to denote a proper transformation acting on $\mathcal{F}$ and a congruence transformation acting on $\mathcal{H}$. Moreover, if the congruence transformation is induced by a reflection, we will also refer to the congruence transformation as a reflection.

For two points $\boldsymbol{u}$ and $\boldsymbol{v}$ in $\mathcal{H}$, define the function $d(\boldsymbol{u}, \boldsymbol{v})$ as

$$d(\boldsymbol{u}, \boldsymbol{v}) = -<q(\boldsymbol{u}), q(\boldsymbol{v})>/(<q(\boldsymbol{u}), p>.<q(\boldsymbol{v}), p>)$$

It is clear that $d(\boldsymbol{u}, \boldsymbol{v}) = d(\boldsymbol{v}, \boldsymbol{u}) \geq 0$ and $d(\boldsymbol{u}, \boldsymbol{v}) = 0$ if and only if $\boldsymbol{u} = \boldsymbol{v}$. We will refer to the function $d(\boldsymbol{u}, \boldsymbol{v})$ as the quasi-distance of $\boldsymbol{u}$ and $\boldsymbol{v}$.

Let $T$ be a congruence transformation that maps $\boldsymbol{u}$ and $\boldsymbol{v}$ to $\boldsymbol{u}'$ and $\boldsymbol{v}'$ respectively. We have $T(q(\boldsymbol{u})) = \lambda q(\boldsymbol{u}')$ and $T(q(\boldsymbol{v})) = \mu q(\boldsymbol{v}')$. It follows readily that $d(\boldsymbol{u}, \mathrm{v}) = d(\boldsymbol{u}', \boldsymbol{v}')$. So the quasi-distance function is invariant under the congruence transformations.

<u>**Proposition 3**</u>: *(a) For any two distinct points* **u** *and* **v** *in* $\mathcal{H}$, *there is a reflection that maps* **u** *to* **v** *and that maps other points on the line through* **u** *and* **v** *to points on the same line.*

*(b) For any two lines L and M in* $\mathcal{H}$, *there is a reflection that maps L to M.*

*(c) Let* **u** *be a point on the line L and* **v** *be a point on the line M, there is a congruence transformation that maps L to M and* **u** *to* **v**.

*(d) Let* $R = [\mathbf{uv}, \infty)$ *and* $S = [\mathbf{xy}, \infty)$ *be two rays in* $\mathcal{H}$. *There is a congruence transformation mapping* **u** *to* **x** *and the ray R to the ray S.*

*Proof.* Let $\lambda$ be the number $> 0$ such that $<q(\boldsymbol{u}), p> = <\lambda q(\boldsymbol{v}), p>$. The cycle $r = q(\boldsymbol{u}) - \lambda q(\boldsymbol{v})$ has norm $-2\lambda<q(\boldsymbol{u}), q(\boldsymbol{v})>$ which is $> 0$. Moreover, $<r, p> = 0$ by the choice of $\lambda$, so $r$ is cycle orthogonal to $p$. The reflection defined by the cycle $r$ is a proper transformation mapping $q(\boldsymbol{u})$ to $\lambda q(\boldsymbol{v})$, so it induces a congruence transformation of order 2 mapping $\boldsymbol{u}$ to $\boldsymbol{v}$. Moreover, if $\ell$ is a cycle corresponding to the line passing through $\boldsymbol{u}$ and $\boldsymbol{v}$, then $\ell$ is orthogonal to both $q(\boldsymbol{u})$ and $q(\boldsymbol{v})$, and therefore $\ell$ is orthogonal to $r$. Accordingly, $\ell$ is invariant under the reflection defined by $r$, and the line defined by $\ell$ is stable under the action of such a reflection. This proves (a).

For (b), let $\ell$ and $m$ be two cycles corresponding to $L$ and $M$, respectively. Because scalar multiples of the same cycle define the same line, we can assume that $\ell$ and $m$ have the same norm by scaling. However, because $\mathrm{norm}(\ell + m) + \mathrm{norm}(\ell - m) = 2(\mathrm{norm}(\ell) + \mathrm{norm}(m)) > 0$, either the cycle $(\ell + m)$ or the cycle $(\ell - m)$ must have strictly positive norm. Since we are free to replace $m$ by $-m$, we can assume that the cycle $(\ell - m)$ has norm $> 0$. The reflection defined by the cycle $(\ell - m)$ maps $\ell$ to $m$, and consequently induces a congruence transformation of order 2 that maps points on the line $L$ to points on the line $M$.

The statement of (c) follows by combining (b) and (a).

From (c), we know there is a congruence transformation that maps $\boldsymbol{u}$ to $\boldsymbol{x}$ and the line $L$ of $\boldsymbol{u}$ and $\boldsymbol{v}$ to the line $M$ of $\boldsymbol{x}$ and $\boldsymbol{y}$. If the transformation maps $\boldsymbol{v}$ to the same side of $\boldsymbol{x}$ as $\boldsymbol{y}$, then we are done. Otherwise, we can compose it with the reflection in (a) to get the desired result. ■

__Proposition 4__: *Let L be a line in H.*

*(a) There is a unique reflection that fixes all points on L and exchanges the two sides of L. Any congruence transformation that fixes all points of L is either the identity transformation or the reflection that exchanges the two sides of L.*

*(b) For any point **u** on L, there is a reflection that fixes **u** and exchanges the two sides of L defined by **u**.*

*(c) If a congruence transformation T maps a line L to itself, and fixes two distinct points on L, then T must fix all points on L. If T fixes just a point **u** on L, then the action of T on the line L must be the same as the reflection in (b).*

*(d) Let **u**, **v**, **w** are three points on a line L such that the segment [**u**, **v**] is congruent to the segment [**u**, **w**]. Then either **v** and **w** are on different sides of **u**, or **v** = **w**.*

*(e) For any two distinct points **v**, **w** on the same side of **u** in a line, we always have d(**u**, **v**) ≠ d(**u**, **w**).*

*Proof.* If a point $\boldsymbol{u}$ is fixed under a congruence transformation $T$, observe that we must have $T(q(\boldsymbol{u})) = q(\boldsymbol{u})$. That is because $\boldsymbol{u}$ is fixed if and only if $T(q(\boldsymbol{u})) = \lambda q(\boldsymbol{u})$ for some positive coefficient $\lambda$, at the same time that $<T(q(\boldsymbol{u})), p> = <q(\boldsymbol{u}), p>$ because $T$ is an orthogonal transformation fixing the cycle $p$.

Let $\ell$ be a cycle corresponding to the line $L$. Any congruence transformation $T$ that fixes the points of $L$ must send $\ell$ to a scalar multiple of itself. Because $T$ is orthogonal, $T(\ell)$ is either $\ell$ or $-\ell$. In the first case, $T$ is the identity transformation. In the second case, $T$ is the reflection defined by $\ell$. That reflection maps $\ell$ to $-\ell$, and leaves invariant any cycle orthogonal to $\ell$. Accordingly, such a reflection fixes all points of $L$ and exchanges the two sides of $L$, and it is uniquely determined by these two properties. That proves (a).

The space generated by $p$, $q(\boldsymbol{u})$ and $\ell$ is a regular 3-dimensional subspace isometric to [1, –1, 1], whose orthogonal complement is generated by a nonisotropic cycle $n$ of positive norm. We can assume by scaling that $n$ has norm 1. The reflection defined by $n$ induces a congruence transformation of order 2 that maps $L$ into $L$ and fixes the point $\boldsymbol{u}$. Moreover, for any other point $\boldsymbol{v}$ on $L$, this transformation maps $\boldsymbol{v}$ to $\boldsymbol{w}$ such that for some positive coefficient $\lambda$, we have $\lambda q(\boldsymbol{w}) = q(\boldsymbol{v}) - 2 <n, q(\boldsymbol{v})>n$. By taking the cycle product with $n$ on both sides, we have $<n, \lambda q(\boldsymbol{w})> = -<n, q(\boldsymbol{v})>$, which means the transformation exchanges the two sides of $\boldsymbol{u}$ in $L$. That proves (b).

Now let $\boldsymbol{u}$ and $\boldsymbol{v}$ be distinct points on $L$ that are fixed by a congruence transformation $T$. $T$ fixes $q(\boldsymbol{u})$, $q(\boldsymbol{v})$ and $p$, so if $\ell$ is a cycle orthogonal to all three of them, then $T(\ell) = \ell$ or $-\ell$. In the first case, $T$ is the identity transformation. In the second case, $T$ is the reflection that fixes all points on $L$ while exchanging the two side of $L$. In either case, all points of $L$ are fixed under the action of $T$.

Assume now that $T$ fixes just a point $\boldsymbol{u}$ on $L$. Because $L$ is stable under the action by $T$, we must have $T(\ell) = \ell$ or $-\ell$. If we only focus on the action of $T$ on $L$, then by composing $T$ with the reflection defined by $\ell$ as needed (such a reflection leaves all points of $L$ invariant), we can assume that $T(\ell) = \ell$. Because $T(q(\boldsymbol{u})) = q(\boldsymbol{u})$ and $T(p) = p$, $T$ must be the reflection defined by a cycle $n$ orthogonal to $p$, $q(\boldsymbol{u})$ and $\ell$. This is the reflection that fixes $\boldsymbol{u}$ and exchanges the two sides of $\boldsymbol{u}$. That proves (c).

Let $T$ be a congruence transformation that maps $[\boldsymbol{u}, \boldsymbol{v}]$ to $[\boldsymbol{u}, \boldsymbol{w}]$. $T$ obviously maps $L$ to itself, and fixes $\boldsymbol{u}$. Hence $T$ must either be the identity, in which case $\boldsymbol{v} = \boldsymbol{w}$, or $T$ flips the two sides of $\boldsymbol{u}$. That proves (d).

To show (e), it is sufficient to consider the case of one specific line, because quasi-distance is invariant under congruence transformations and congruence transformations act transitively on the set of lines.

Consider the line defined by the 1-cycle $\boldsymbol{j}.\boldsymbol{X}$, where $\boldsymbol{j}$ is a nonzero vector in $\boldsymbol{E}$ such that $\boldsymbol{i}.\boldsymbol{j} = 0$. This line consists of all points in $E$ of the form $-2t\boldsymbol{i}$ where $t$ runs through all numbers $> 0$. On this line, let the points $\boldsymbol{u}$, $\boldsymbol{v}$ and $\boldsymbol{w}$ correspond to $-2a\boldsymbol{i}$, $-2b\boldsymbol{i}$ and $-2c\boldsymbol{i}$. To say that $\boldsymbol{v}$ and $\boldsymbol{w}$ are on the same side of $\boldsymbol{u}$ is the same as saying that the numbers $b$ and $c$ are both larger or both smaller than $a$.

We have $d(\boldsymbol{u}, \boldsymbol{v}) = 4(a - b)^2/4ab = (a - b)^2/ab$. Similarly, $d(\boldsymbol{u}, \boldsymbol{w}) = (a - c)^2/ac$. Accordingly, $d(\boldsymbol{u}, \boldsymbol{v}) = d(\boldsymbol{u}, \boldsymbol{w})$ if and only if $c(a - b)^2 = b(a - c)^2$, or after multiplying out and rearranging terms, if and only if $(a^2 - bc)(b - c) = 0$. Because $b$ and $c$ are distinct numbers either both larger than $a$ or both smaller than $a$, this equation is impossible. That proves statement (e). ■

From the above properties, we can readily define the concept of congruence and verify the standard Hilbert axioms for congruence. Given two points $\boldsymbol{u}$ and $\boldsymbol{v}$ in $\mathcal{H}$, the segment $[\boldsymbol{u}, \boldsymbol{v}]$ is the subset of $\mathcal{H}$ consisting of $\boldsymbol{u}$, $\boldsymbol{v}$ and all points between $\boldsymbol{u}$ and $\boldsymbol{v}$. We say the segment $[\boldsymbol{u}, \boldsymbol{v}]$ is congruent to the segment $[\boldsymbol{u}', \boldsymbol{v}']$ if there is a congruence transformation mapping $\boldsymbol{u}$ to $\boldsymbol{u}'$ and $\boldsymbol{v}$ to $\boldsymbol{v}'$. Such a transformation will of course map all points between $\boldsymbol{u}$ and $\boldsymbol{v}$ to all points between $\boldsymbol{u}'$ and $\boldsymbol{v}'$ because congruence transformations preserve the betweenness relationship. Based on the foregoing discussion, it follows readily that two segments are congruent if and only if they have the same quasi-distance.

By an angle we understand an ordered pair of two rays $[\boldsymbol{uv}, \infty)$ and $[\boldsymbol{uw}, \infty)$ issuing from the same point $\boldsymbol{u}$. We say two angles are congruent if there is a congruence transformation mapping the two rays of one angle to the two rays of the other angle (in the same order).

We can now proceed to the axioms of congruence.

*(Congruence Axiom 1) A segment is congruent to itself. If two segments are congruent to a third one, they are congruent to each other.*

*(Congruence Axiom 2) An angle is congruent to itself. If two angles are congruent to a third one, they are congruent to each other.*

Under our definition, congruence of segments and angles is obviously an equivalence relation. ■

*(Congruence Axiom 3) Let **u** and **v** be two distinct points in $\mathcal{H}$. If **w** is any point in and L is any line passing through **w**, there is on each side of **w** in L a unique point **t** such that the segment [**w, t**] is congruent to the segment [**u**, **v**].*

From Proposition 3, we know that there is a congruence transformation mapping ***u*** to ***w*** and the line through ***u*** and ***v*** to the line *L*. Such a transformation will map the point ***v*** to a point ***t*** on one of the two sides of ***w*** in *L*. The side where ***t*** falls will completely determine such a transformation, for if we have two different transformations that map ***v*** to different points on the same side of ***w***, then by composing the inverse of one transformation with the other, we have a congruence transformation that maps *L* to *L* and fixes ***u***. According to Proposition 4, such a transformation must either be the identity or flip the two sides of ***u*** in *L*. But we just construct a transformation that is not the identity and does not flip the two sides, a contradiction.

Whatever side of ***w*** the point ***t*** falls on, we can compose the original congruence transformation with the reflection that fixes ***w*** and flips the two sides of ***w*** in *L*, giving us another point on the other side with similar property. ■

*(Congruence Axiom 4) On a line L let [**u**, **v**] and [**v**, **w**] be two segments with only the point **v** in common. Suppose on a line M we have segments [**x**, **y**] and [**y**, **z**] with only the point **y** in common and that are congruent to [**u**, **v**] and [**v**, **w**] respectively. In that case, the segment [**u**, **w**] is congruent to the segment [**x**, **z**].*

Let $T$ be a congruence transformation that maps the line $L$ to the line $M$ and sends the point $\boldsymbol{v}$ to the point $\boldsymbol{y}$. The hypothesis implies that $\boldsymbol{u}$ and $\boldsymbol{w}$ are on different sides of $\boldsymbol{v}$ in $L$, and $\boldsymbol{x}$ and $\boldsymbol{z}$ are on different sides of $\boldsymbol{y}$ in $M$. By composing $T$ with an appropriate reflection, we can assume that $T$ maps $\boldsymbol{u}$ to $\boldsymbol{x}$. In that case, $T$ must map $\boldsymbol{w}$ to $\boldsymbol{z}$, making the segment $[\boldsymbol{u}, \boldsymbol{w}]$ congruent to the segment $[\boldsymbol{x}, \boldsymbol{z}]$. ■

*(Congruence Axiom 5) Let α be an angle defined by two rays R and S. Let T be a ray that is part of a line L. There is on each side of the line L a unique ray making an angle with T that is congruent to the angle α.*

The existence of such a ray follows from Proposition 3(d) and Proposition 4(a). Uniqueness also follows from Proposition 4(a), namely a congruence transformation that fixes all the points of a line $L$ and does not exchange the two sides of $L$ in $\mathcal{H}$ must be the identity transformation. ■

*(Congruence Axiom 6 – The Side-Angle-Side Congruence Axiom) Let* $\mathbf{u}, \mathbf{v}, \mathbf{w}$ *and* $\mathbf{x}, \mathbf{y}, \mathbf{z}$ *be two sets of three non-collinear points in* $\mathcal{H}$*. If the segments* $[\mathbf{u}, \mathbf{v}]$ *and* $[\mathbf{u}, \mathbf{w}]$ *are congruent to the segments* $[\mathbf{x}, \mathbf{y}]$ *and* $[\mathbf{x}, \mathbf{z}]$ *respectively, and the angle determined by the rays* $[\mathbf{uv}, \infty)$ *and* $[\mathbf{uw}, \infty)$ *is congruent to the angle determined by the rays* $[\mathbf{xy}, \infty)$ *and* $[\mathbf{xz}, \infty)$*, then there is a congruence transformation mapping* $\mathbf{u}, \mathbf{v}, \mathbf{w}$ *to* $\mathbf{x}, \mathbf{y}, \mathbf{z}$ *respectively.*

Let $T$ be a congruence transformation that maps the segment $[\boldsymbol{u}, \boldsymbol{v}]$ to the segment $[\boldsymbol{x}, \boldsymbol{y}]$. $T$ obviously maps the line $L$ through $\boldsymbol{u}$ and $\boldsymbol{v}$ to the line $M$ through $\boldsymbol{x}$ and $\boldsymbol{y}$. The congruence transformation $T$ must also map the ray $[\boldsymbol{uw}, \infty)$ to one of the two rays on each side of $M$ that, together with the ray $[\boldsymbol{xy}, \infty)$, makes an angle congruent to the angle determined by the rays $[\boldsymbol{uv}, \infty)$ *and* $[\boldsymbol{uw}, \infty)$. By composing $T$ with the reflection that fixes the line $M$ and exchanges the two sides of $M$, we can assume that $T$ maps the ray $[\boldsymbol{uw}, \infty)$ to the same side of $M$ as the point $\boldsymbol{z}$. In that case, it is clear that $T$ maps the ray $[\boldsymbol{uw}, \infty)$ to the ray $[\boldsymbol{xz}, \infty)$. By construction, the segment $[\boldsymbol{x}, \boldsymbol{z}]$ is congruent to $[\boldsymbol{x}, T(\boldsymbol{w})]$, and $\boldsymbol{z}$ and $T(\boldsymbol{w})$ are on the same side of $\boldsymbol{x}$. By Proposition 4(d), we must have $T(\boldsymbol{w}) = \boldsymbol{z}$. ■

**6. THE HYPERBOLIC PLANE – PARALLEL AXIOM.** We can now consider the parallel axiom for our hyperbolic plane.

*(Parallel Axiom – Hyperbolic Version) Let L be a line in* $\mathcal{H}$ *and let* $\boldsymbol{u}$ *be a point not on L. There are two or more lines passing through* $\boldsymbol{u}$ *that do not intersect with L.*

We will show that there are in fact infinitely many lines passing through the point $\boldsymbol{u}$ that do not intersect with the line $L$. Lines that do not intersect are also said to be parallel.

In general, let $\ell$ and $\mathscr{m}$ be nonisotropic cycles that define two distinct lines in $\mathcal{H}$. We claim that the lines defined by $\ell$ and $\mathscr{m}$ do not intersect if and only if the 2-dimensional cycle subspace generated by $\ell$ and $\mathscr{m}$ is isotropic.

Indeed, we will show that the lines defined by $\ell$ and $\mathscr{m}$ intersect if and only if the space generated by $\ell$ and $\mathscr{m}$ is anisotropic, or equivalently, that its orthogonal complement is an Artinian plane.

It is clear that if these lines intersect at a point $\boldsymbol{w}$, then $q(\boldsymbol{w})$ is cycle orthogonal to both $\ell$ and $\mathscr{m}$, and hence belongs to their orthogonal complement. That orthogonal complement has dimension 2 and also includes $p$. The orthogonal complement is therefore generated by $p$ and $q(\boldsymbol{w})$. It is regular and isotropic, and hence is an Artinian plane.

Now consider the converse. We want to show that if the orthogonal complement to $\ell$ and $\mathscr{m}$ is an Artinian plane, the two lines must intersect in $\mathcal{H}$. Such an Artinian plane must contain two isotropic lines not orthogonal to $p$. These isotropic lines therefore are multiples of the normalized 2-cycles $q(\boldsymbol{v})$ and $q(\boldsymbol{w})$ for two points $\boldsymbol{v}$ and $\boldsymbol{w}$ in $E$. It is straight-forward to verify that we must have $q(\boldsymbol{w}) = q(\boldsymbol{v}) - 2 \langle p, q(\boldsymbol{v})\rangle\, p$, so $\langle p, q(\boldsymbol{v})\rangle$ and $\langle p, q(\boldsymbol{w})\rangle$ must have different signs, forcing one of the points $\boldsymbol{v}$ and $\boldsymbol{w}$ to be in $\mathcal{H}$. That point is the intersection of the two lines defined by $\ell$ and $\mathscr{m}$. Our claim is now proved.

Let $D$ be the orthogonal complement of $p$ in the space $\mathcal{F}$ of all cycles. $D$ is a space of dimension 3, isometric to the quadratic space $[1, 1, -1]$. The cycles in $D$ orthogonal to $q(\boldsymbol{u})$ form a positive definite subspace of dimension 2, because those cycles are the orthogonal complement in $\mathcal{F}$ of the Artinian plane generated by $p$ and $q(\boldsymbol{u})$. The nonzero vectors in that subspace define lines that pass through the point $\boldsymbol{u}$, with scalar multiples of the same vector representing the same line.

Now look at the projective plane $\boldsymbol{P}(D)$ associated to $D$. In that projective plane, the lines passing through the point $\boldsymbol{u}$ correspond to points on the projective line $\mathcal{U}$ defined by the 2-dimensisonal subspace of cycles orthogonal to $q(\boldsymbol{u})$. Let $\ell$ be the cycle that defines the given line $L$ in $\mathcal{H}$, and let the point $\Lambda$ be the image of $\ell$ in $\boldsymbol{P}(D)$. If a cycle $\mathcal{m}$ defines a line passing through $\boldsymbol{u}$, its image in $\boldsymbol{P}(D)$ is a point $\mu$ on the line $\mathcal{U}$, and the projective line joining $\Lambda$ and $\mu$ represents the 2-dimensional subspace generated by $\ell$ and $\mathcal{m}$. Such a subspace is isotropic if and only if the line joining $\Lambda$ and $\mu$ intersects the projective conic defined by the cycle pairing in $D$. That cycle pairing conic is isomorphic to the projective conic $X^2 + Y^2 - Z^2 = 0$, and therefore has an infinite number of points.

For each point on the conic, consider the line joining $\Lambda$ and that point. That line will intersect the line $\mathcal{U}$ in a point because any two lines in a projective plane intersect. In light of our foregoing discussion, that point represents a line through $\boldsymbol{u}$ that is parallel to $L$.

Because the point $\Lambda$ is not on the conic, the lines joining $\Lambda$ with points on $\mathcal{U}$ and intersecting with the cycle pairing conic would correspond roughly to half of that conic. Except for the two tangents to the conic drawn from $\Lambda$, each line through $\Lambda$ that intersects with the conic will do so in exactly two distinct points. (In terms of bilinear algebra, such a line represents an Artinian plane, and therefore has exactly two isotropic points.) Accordingly, there are an infinite number of lines through $\boldsymbol{u}$ that do not intersect with $L$. ■

______________________________________________________________________.